\documentclass[12pt]{article}
\usepackage{amsfonts, amsmath, amssymb,amsthm,bbm,comment,fancyhdr,float,makeidx,mathrsfs,verbatim}
\usepackage{latexcad}

\normalsize

\newcommand{\seq}{\subseteq}

\newcommand{\vs}{\vspace*}

\newcommand{\nin}{\noindent}

\newtheorem{mthm}{Theorem}[section]
\newtheorem{mylem}[mthm]{Lemma}
\newtheorem{myprn}[mthm]{Proposition}
\newtheorem{mycor}[mthm]{Corollary}
\newtheorem{mydef}[mthm]{Definition}
\newtheorem{myrem}[mthm]{Remark}
\newtheorem{mycon}[mthm]{Construction}
\newtheorem{myeg} [mthm]{Example}
\newtheorem{myque} [mthm]{Question}
\newtheorem{myalg} [mthm]{Algorithm}
\newenvironment{thm}{\begin{mthm}}{\end{mthm}}
\newenvironment{lem}{\begin{mylem}}{\end{mylem}}
\newenvironment{cor}{\begin{mycor}}{\end{mycor}}
\newenvironment{prop}{\begin{myprn}}{\end{myprn}}
\newenvironment{mdef}{\begin{mydef}\rm}{\end{mydef}}
\newenvironment{rem}{\begin{myrem}\rm}{\end{myrem}}

\newenvironment{ex}{\begin{myeg}\rm}{\end{myeg}}
\newenvironment{prof}{\noindent $Proof.$ \rm}{\hfill $\Box$}

\newenvironment{alg}{\begin{myalg}\rm}{\end{myalg}}

\def \nin {\noindent}

\def \Lemma #1 {\vs{3mm}\nin {\bf Lemma #1} \it}
\def \Prop #1 {\vs{3mm}\nin {\bf Proposition #1} \it}
\def \Th #1 {\vs{3mm}\nin {\bf Theorem #1} \it}
\def \Cor #1 {\vs{3mm}\nin {\bf Corollary #1} \it}
\def \Ex #1 {\vs{3mm}\nin {\bf Example #1} \it}

\def \part #1 {\hfil\break\hglue 12pt {\rm (#1)~}}

\def\fs{\footnotesize}

\title{
\bf\LARGE  On the $(n, d)^{th}$ $f$-Ideals \thanks{This research is supported by the National Natural
Science Foundation of China (Grant No. 11271250). } }
\author{Jin Guo\thanks{Corresponding author. guojinecho@163.com}, \,\, {Tongsuo Wu\thanks{ tswu@sjtu.edu.cn}}\\
 {\small Department of Mathematics, Shanghai Jiaotong University}
}

\date{}
\begin{document}
\baselineskip=16pt \maketitle

\begin{center}
\begin{minipage}{11cm}

 \vs{3mm}\nin{\small\bf Abstract.} {\fs A square-free monomial ideal $I$ is called an {\it $f$-ideal}, if both $\delta_{\mathcal{F}}(I)$ and $\delta_{\mathcal{N}}(I)$ have the same $f$-vector, where $\delta_{\mathcal{F}}(I)$ ($\delta_{\mathcal{N}}(I)$, respectively) is the facet (Stanley-Reisner, respectively) complex related to $I$. In this paper, we introduce the concepts of perfect set containing $k$ and perfect set without $k$. We study the $(n, d)^{th}$ perfect sets and show that $V(n, d) \neq \emptyset$ for $d \geq 2$ and $n \geq d+2$. Then we give some algorithms to construct $(n, d)^{th}$ $f$-ideals and show an upper bound of the $(n, d)^{th}$ perfect number.
}

\vs{3mm}\nin {\small\bf Key Words and phrases:} {\small  perfect set containing $k$; perfect set without $k$; $f$-ideal; unmixed $f$-ideal; perfect number}

\vs{3mm}\nin {\small\bf 2010 Mathematics Subject Classification: 13P10, 13F20, 13C14, 05A18. }

\end{minipage}
\end{center}

\vs{4mm} \section{Introduction }

Throughout the paper, for a set $A$, we use $A_d$ to denote the set  of  the subsets of $A$ with cardinality $d$. For a monomial ideal $I$ of $S$, let $sm(I)$ be the set of square-free monomials in $I$.
As we know, there is a natural bijection between $sm(S)$ and $2^{[n]}$, denoted by
$$\sigma: \, x_{i_1}x_{i_2}\cdots x_{i_k} \mapsto \{i_1, i_2, \ldots, i_k\},$$
where $[n]=\{1,2,\ldots,n\}$ for a positive integer $n$.
A square-free monomial $u$ is called {\it covered} by a square-free monomial $v$, if $u \,|\, v$ holds.
For other concepts and notations, see references \cite{AM, Eisenbud, GWL, HH, Villareal,Zariski}.

Constructing free resolutions of a monomial ideal is one of the core problems in commutative algebra. A main approach to the problem is by taking advantage of the properties of a simplicial complex, so it is important to have a research on the properties of the complex corresponding to the related ideals,  see for example, references \cite{CF, Faridi, HHZ, Zheng}. There is an important class of ideals called $f$-ideals,  whose facet complex $\delta_{\mathcal{F}}(I)$  and Stanley-Reisner complex $\delta_{\mathcal{N}}(I)$ have the same $f$-vector, where $\delta_{\mathcal{F}}(I)$ is generated by the set $\sigma(G(I))$, and $\delta_{\mathcal{N}}(I) = \{ \sigma(g) \,|\, g \in sm(S) \setminus sm(I)\}$.  Note that the $f$-vector of a complex $\delta_{\mathcal{N}}(I)$, which is not easy to calculate in general,  is essential in the computation of  the Hilbert series of $S/I$. Since the correspondence of the complex $\delta_{\mathcal{F}}(I)$ and the ideal $I$ is direct and clear, it is more easier to calculate the $f$-vector of $\delta_{\mathcal{F}}(I)$. So, it is convenient to calculate the Hilbert series and study other corresponding properties of $S/I$ while $I$ is an $f$-ideal.

The formal definition of an $f$-ideal first appeared in \cite{AAAB}, and  it was then studied in \cite{AMBZ}. In \cite{GWL}, the authors characterized the $f$-ideals of degree $d$, as well as the $f$-ideals in general case. They introduced a bijection between square-free monomial ideals of degree $2$ and simple graphs, and showed that $V(n, 2) \neq \emptyset$ for any $n \geq 4$, where $V(n, d)$ is the set of $f$-ideals of degree $d$ in $S=K[x_1, \ldots, x_n]$. The structure of $V(n, 2)$ was determined, and the characterization of the unmixed $f$-ideals is also studied in \cite{GWL}.

In this paper, we give another characterization of unmixed $f$-ideals in part two. In part three, we generalize the aforementioned result of \cite{GWL} by showing that $V(n, d) \neq \emptyset$ for general $d \geq 3$ and $n \geq d+2$. In part four, we introduce some algorithms to construct  $(n, d)^{th}$ $f$-ideals, and we show an upper bound of the $(n, d)^{th}$ perfect number in part five. In part six, we show some examples of nonhomogeneous $f$-ideals, which is still open in \cite{GWL}.


The following propositions are needed in this paper.

\begin{prop}\label{homogeneous $f$-Ideal}(\cite{GWL} Theorem 2.4) Let $S=K[x_1, \ldots, x_n]$, and let $I$ be a square-free monomial ideal of $S$ of degree $d$ with the minimal generating set $G(I)$. Then $I$ is an $f$-ideal if and only if $\, G(I)$ is $(n, d)^{th}$ perfect  and $|G(I)| = \frac{1}{2}C_n^d$ holds true.
\end{prop}

\begin{prop}\label{V(n,2) nonempty}(\cite{GWL} Proposition 3.3)  $V(n,2) \neq \emptyset$ if and only if $n=4k$ or $n=4k+1$ for some positive integer $k$.
\end{prop}

\begin{prop}\label{unmixed homogeneous $f$-Ideal}(\cite{GWL} Proposition 5.3) Let $S=K[x_1, \ldots, x_n]$. If $I$ is an $f$-ideal of $S$ of degree $d$, then $I$ is unmixed if and only if $sm(S)_d \setminus G(I)$ is lower perfect.
\end{prop}

\vs{3mm}In \cite{GWL}, a method for finding an $(n, 2)^{th}$ perfect set with the smallest cardinality is provided, namely, first, decompose the set $[n]$ into a disjoint union of two subsets $B$ and $C$ {\it uniformly}, i.e.,  $||B|-|C|| \leq 1$ holds; then set $A= \{x_ix_j \,|\, i, j \in B,$ or $i, j \in C\}$. Finally,  $A$ is an $(n, 2)^{th}$ perfect set whose cardinality is equal to the $(n, 2)^{th}$ perfect number $N_{(n, 2)}$, where
\begin{equation}\label{eq: perfect number new understanding}
N_{(n,2)} = \begin{cases}
k^2-k,  &  \text{if } n=2k; \\
k^2,  &  \text{if } n=2k+1.
\end{cases}
\end{equation}
Note that any set $D$ with $A \seq D \seq sm(S)_2$ is also an $(n,2)^{th}$  perfect set.

\section{ $(n, d)^{th}$ unmixed $f$-ideals }

For a positive integer $d$ greater than $2$, an $(n, d)^{th}$ $f$-ideal may be not unmixed, see Example 5.1 of \cite{GWL} for a counterexmple. So, it is interesting to characterize the unmixed $f$-ideals. In this section, we show a characterization of unmixed $f$-ideals by the corresponding simplicial complex, by taking advantage of the bijection $\sigma$ between square-free monomial ideals and  simplicial complexes.

Recall that a simplicial complex is a {\it $d$-flag complex} if all of its minimal non-faces contain $d$ elements.  Recall that $\Delta^{\vee}$ denotes the Alexander dual of a simplicial complex $\Delta$, see \cite{HH} for details.

\begin{prop}\label{unmixed f-ideal} Let $S=K[x_1, \ldots, x_n]$, and let $I$ be a square-free monomial ideal of $S$ of degree $d$. $I$ is an $(n, d)^{th}$ unmixed $f$-ideal if and only if the followings hold:\\
$(1)$ $|G(I)|= C_n^d/2$;\\
$(2)$ $dim\,\delta_{\mathcal{F}}(I)^{\vee}=n-d-1$;\\
$(3)$ $\langle \sigma(u) \,|\, u \in sm(S)_d \setminus G(I) \rangle$ is a $d$-flag complex.
\end{prop}

\begin{prof}We claim that the following two results hold true: First, the condition (2) holds if and only if $G(I)$ is lower perfect. Second, the condition (3) holds if and only if $G(I)$ is upper perfect and $sm(S)_d \setminus G(I)$ is lower perfect. If the above two results hold true, then it is easy to see that the conclusion holds by Proposition \ref{homogeneous $f$-Ideal} and Proposition \ref{unmixed homogeneous $f$-Ideal}.

For the first claim, if $G(I)$ is lower perfect, then for each minimal non-face $F$ of $\delta_{\mathcal{F}}(I)$, $|F| \geq d$ holds. By the definition of the Alexander dual, $G$ is a face of $\delta_{\mathcal{F}}(I)^{\vee}$ if and only if $[n] \setminus G$ is a non-face of $\delta_{\mathcal{F}}(I)$. So, for each facet $L$ of $\delta_{\mathcal{F}}(I)^{\vee}$,  $|L| \leq n-d$. Since $|G(I)| \neq C_n^d$, there exists some non-face of $\delta_{\mathcal{F}}(I)$ with cardinality $d$, so there exists some facet of $\delta_{\mathcal{F}}(I)^{\vee}$ with cardinality $n-d$. Thus $dim(\delta_{\mathcal{F}}(I)^{\vee})=n-d-1$.

Conversely, assume $dim(\delta_{\mathcal{F}}(I)^{\vee})=n-d-1$. By a similar argument, one can see that the smallest cardinality of non-faces of $\delta_{\mathcal{F}}(I)$ is $d$, hence $G(I)$ is lower perfect.

For the second claim, if $sm(S)_d \setminus G(I)$ is lower perfect, then for the complex $\Delta= \langle \sigma(u) \,|\, u \in sm(S)_d \setminus G(I) \rangle$, the cardinality of a non-face is not less than $d$. Since $G(I)$ is upper perfect,  for each non-face $F$ of $\Delta$, there exists $v \in G(I)$ such that $\sigma(v) \seq F$. Note that $\sigma(v)$ is a non-face of $\Delta$, so all the minimal non-faces of $\Delta$ have cardinality $d$. Hence $\Delta$ is a $d$-flag complex.

Conversely,
assume that $\Delta = \langle \sigma(u) \,|\, u \in sm(S)_d \setminus G(I) \rangle$ is a $d$-flag complex. In a similar way, one can see that $G(I)$ is upper perfect and $sm(S)_d \setminus G(I)$ is lower perfect.
\end{prof}

\section{ Existence of $(n, d)^{th}$ $f$-ideals }

 For a subset $M$ of $sm(S)_d$, denote $M^{'}=\{\sigma^{-1}(A) \,|\, A=[n] \setminus \sigma(u)$ for some $u \in M\}$.
The following lemma is essential in the proof of our main result in this section.

\begin{lem}\label{complement lemma} $M$ is a perfect subset of $sm(S)_d$ if and only if $M^{'}$ is a perfect subset of $sm(S)_{n-d}$.
\end{lem}

\begin{prof}For the necessary part, if $M$ is a subset of $sm(S)_d$, then it follows from definition that $M^{'}$ is a subset of $sm(S)_{n-d}$. In order to check that $M^{'}$ is upper perfect, we will show for each monomial $u \in sm(S)_{n-d+1}$ that $u \in \sqcup(M^{'})$ holds. This is equivalent to showing that there exists some $v\in M^{'}$, such that $\sigma(v) \seq \sigma(u)$ holds. In fact, since $M$ is lower perfect, for the monomial $u^{'}=\sigma^{-1}([n] \setminus \sigma(u)) \in sm(S)_{d-1}$, there exists some $w \in M$ such that $u^{'} \,|\, w$ holds. Hence $\sigma(u^{'}) \seq \sigma(w)$ holds. Now let $v= \sigma^{-1}([n] \setminus \sigma(w))$ and then it is easy to see that $\sigma(v)=[n] \setminus \sigma(w) \seq [n] \setminus \sigma(u^{'}) = \sigma(u)$ hold. This shows that $M^{'}$ is upper perfect. In a similar way, one can prove that $M^{'}$ is lower perfect.

The sufficient part is similar to prove, and we omit the details.
\end{prof}

\vs{3mm}By the proof of the above lemma, one can see that $M$ is an upper (lower, respectively) perfect subset of $sm(S)_d$ if and only if $M^{'}$ is a lower (upper, respectively) perfect subset of $sm(S)_{n-d}$.

\begin{cor}\label{complement cor}If $I$ is a square-free monomial ideal of $S$ of degree $d$, then $I$ is an $(n, d)^{th}$ $f$-ideal if and only if $|G(I)|=C_n^d/2$ and $G(I)^{'}$ is a perfect subset of $sm(S)_{n-d}$.
\end{cor}

\vs{3mm}Denote $sm(S\{\check{k}\})_d=\{u\in sm(S)_d \,|\; x_k \nmid u\}$, and $sm(S\{k\})_d=\{u\in sm(S)_d \,|\; x_k | u\}$. For a subset $X = \{i_1, \ldots, i_j\}$ of $[n]$, denote $$sm(S\{\check{X}\})_d=\{u\in sm(S)_d \,|\; x_k \nmid u\,\,for\, every\,\, k \in X\},$$ and let $sm(S\{X\})_d=\{u\in sm(S)_d \,|\; x_k \mid u$ for every $k \in X\}$.

\begin{mdef}\label{perfect containing without k} For a subset $M$ of $sm(S\{\check{k}\})_d$,  if $sm(S\{\check{k}\})_{d+1} \seq \sqcup(M)$ holds, then $M$ is called {\it upper perfect without $k$}. Dually, a subset $M$ of $sm(S\{\check{k}\})_d$ is called {\it lower perfect without $k$}, if $sm(S\{\check{k}\})_{d-1} \seq \sqcap(M)$ holds. A subset $M$ of $sm(S\{k\})_d$ is called {\it upper perfect containing $k$}, if $sm(S\{k\})_{d+1} \seq \sqcup(M)$ holds;  a subset $M$ of $sm(S\{k\})_d$ is called {\it lower perfect containing $k$}, if $sm(S\{k\})_{d-1} \seq \sqcap(M)$ holds. If $M$ is not only upper but also lower perfect without $k$, then $M$ is called {\it perfect without $k$}. Similarly, if $M$ is both upper and lower perfect containing $k$, then $M$ is called {\it perfect containing $k$}.
\end{mdef}

\vs{3mm}For a subset $X$ of $[n]$, we can define the upper perfect (lower perfect, perfect, respectively) set without $X$ (containing $X$) similarly. For a subset $A$ of $sm(S)_d$, let $A\{\check{X}\} = A \cap sm(S\{\check{X}\})_d$, and let $A\{X\} = A \cap sm(S\{X\})_d$.

\begin{prop}\label{upper without lower containing} Let $A$  be a subset of $sm(S)_d$, and let $X= \{i_1, \ldots, i_j\}$ be a subset of $[n]$.  Then the  following statements hold:

$(1)$ $A\{\check{X}\} = A\{\check{i_1}\}\{\check{i_2}\} \ldots \{\check{i_j}\}$, and $A\{X\} = A\{i_1\}\{i_2\}\ldots\{i_j\}$;

$(2)$ If $A$ is upper perfect, then $A\{\check{X}\}$ is upper perfect without $X$;

$(3)$ If $A$ is lower perfect, then $A\{X\}$ is lower perfect containing $X$;

$(4)$ If $A$ is upper (lower, respectively) perfect without $X$, then $A^{'}$ is lower (upper, respectively) perfect containing $X$. Furthermore, the converse also holds true.
\end{prop}


\begin{prof}(1) and (2) are easy to see by the corresponding definitions.

In order to prove (3), it is sufficient to show that $A\{k\}$ is a lower perfect set containing $k$ for each $k\in [n]$.
In fact, since $A$ is lower perfect, for each monomial $u \in sm(S\{k\})_{d-1}$,  there exists a monomial $v$ in $A$ such that $u \,|\, v$. Note that $x_k \,|\, u$ holds, so $x_k \,|\, v$ also holds, which implies that $v \in sm(S\{k\})_d$ holds. Hence $A\{k\}$ is a lower perfect set containing $k$.

For (4), we only show that $A^{'}$ is lower perfect containing $k$ when $A$ is upper perfect without $k$, and the remaining implications are similar to prove. In fact, for each monomial $u \in sm(S\{k\})_{n-d-1} \seq sm(S)_{n-d-1}$, $u^{'} \in sm(S)_{d+1}$, note that $x_k \,|\, u$ implies $x_k \nmid u^{'}$ holds true, hence $u^{'} \in sm(S\{\check{k}\})_{d+1}$ also hold. Since $A$ is upper perfect without $k$, there exists a monomial $v \in A$ such that $v \,|\, u^{'}$ holds, hence $u \,|\, v^{'}$ holds, where $v^{'} \in A^{'}$. This completes the proof.
\end{prof}

\begin{rem}\label{not lower without, not upper containing} For a perfect subset $A$ of $sm(S)_d$, $A\{\check{X}\}$ needs not to be a lower perfect set without $X$, and $A\{X\}$ needs not to be an upper perfect set containing $X$, see the following for counter-examples:
\end{rem}

\begin{ex}\label{not upper containing example}Let $S= K[x_1, \ldots, x_6]$,  let $$A = \{x_1x_2x_3, x_1x_2x_4, x_1x_2x_5, x_3x_4x_5, x_1x_2x_6, x_1x_3x_6, x_2x_3x_6, x_4x_5x_6\},$$ and let $B=A \setminus \{x_1x_2x_6\}$. It is easy to see
$$A\{\check{6}\} = B\{\check{6}\} = \{x_1x_2x_3, x_1x_2x_4, x_1x_2x_5, x_3x_4x_5\},$$ $A\{6\} = \{ x_1x_2x_6, x_1x_3x_6, x_2x_3x_6, x_4x_5x_6\}$, and $B\{6\}=\{x_1x_3x_6, x_2x_3x_6, x_4x_5x_6\}$. Also, it is direct to check that both $A$ and $B$ are perfect sets, and that both $A\{\check{6}\}$ and $B\{\check{6}\}$ are perfect sets without $6$. Note that $A\{6\}$ is a perfect set containing $6$, but $B\{6\}$ is not upper perfect.
\end{ex}

\vs{3mm}By Proposition \ref{upper without lower containing}, we have the following example by mapping $A$,  $B$ to $A^{'}$, $B^{'}$ respectively.

\begin{ex}\label{not lower without example}Let $S= K[x_1, \ldots, x_6]$, and let $$A^{'} = \{x_1x_2x_3, x_1x_4x_5, x_2x_4x_5, x_3x_4x_5, x_1x_2x_6, x_3x_4x_6, x_3x_5x_6, x_4x_5x_6\},$$ and $B^{'}=A^{'} \setminus \{x_3x_4x_5\}$. It is easy to see that $$A^{'}\{\check{6}\} = \{x_1x_2x_3, x_1x_4x_5, x_2x_4x_5, x_3x_4x_5\},\,B^{'}\{\check{6}\} = \{x_1x_2x_3, x_1x_4x_5, x_2x_4x_5\},$$ and $A^{'}\{6\} = B^{'}\{6\} = \{x_1x_2x_6, x_3x_4x_6, x_3x_5x_6, x_4x_5x_6\}$. It is direct to check that both $A^{'}$ and $B^{'}$ are perfect sets, and that both $A^{'}\{6\}$ and $A^{'}\{6\}$ are perfect sets containing $6$. Note that $A^{'}\{\check{6}\}$ is a perfect set without $6$, but $B^{'}\{\check{6}\}$ is not lower perfect.
\end{ex}
In order to obtain the main result of this section, we need a further fact and we omit the verification.

\begin{lem}\label{without add containing} Let $S=K[x_1, \ldots, x_n]$, and let $A$ be a subset of $sm(S)_d$. If $A\{\check{k}\}$ is a perfect subset of $sm(S\{\check{k}\})_d$ without $k$, and $A\{k\}$ is a perfect subset of $sm(S\{k\})_d$ containing $k$  for some $k\in [n]$, then $A$ is a perfect subset of $sm(S)_d$.
\end{lem}

\begin{thm}\label{existence of (n, d)^{th} perfect sets}For any integer $d \geq 2$ and any integer $n \geq d+2$, there exists an $(n, d)^{th}$ perfect set with cardinality less than or equal to $C_n^d/2$.

\end{thm}

\begin{prof} We prove the result by induction on $d$.

If $d=2$, the conclusion holds true for any integer $n \geq 4$ by Proposition \ref{V(n,2) nonempty}. In the following, assume $d>2$.

Assume that the conclusion holds true for any integer less than $d$. For $d$, we claim that the conclusion holds true for any integer $n \geq d+2$. We will show the result by induction on $n$.

If $n=d+2$, then $C_n^d = C_n^2$. Note that for any integer $n \geq 4$, there exists an $(n, 2)^{th}$ perfect set $M$, such that $|M| \leq C_n^2/2$. By Lemma \ref{complement lemma}, $M^{'}$ is an $(n, d)^{th}$ perfect set. Note that  $|M^{'}|=|M| \leq C_n^2/2=C_n^d/2$ holds.

Now assume that the conclusion holds true for any integer less than $n$. Then by Lemma 3.8,  it will suffice to show that there is a perfect subset $A$ of $sm(S\{\check{n}\})_d$ without $n$ and a perfect subset $B$ of $sm(S\{n\})_d$ containing $n$, such that $|A| \leq |sm(S\{\check{n}\})_d|/2=C_{n-1}^d/2$ and $|B| \leq |sm(S\{n\})_d|/2=C_{n-1}^{d-1}/2$ hold.

Let $L=K[x_1, \ldots, x_{n-1}]$. Then clearly, $sm(S\{\check{n}\})_d = sm(L)_d$ holds. By induction on $n$, there exists an $(n-1, d)^{th}$ perfect subset $A$ of $sm(L)_d$, such that $|A| \leq C_{n-1}^d/2$. It is easy to see that $A$ is a perfect subset of $sm(S\{\check{n}\})_d$ without $n$.
By induction on $d$, there exists an $(n-1, d-1)^{th}$ perfect subset $B_1$ of $sm(L)_{d-1}$, such that $|B_1| \leq C_{n-1}^{d-1}/2$ holds. Let $B=\{\sigma^{-1}(D) \,|\, D=\sigma(u) \cup \{n\}$ for some $u \in B_1\}$. It is easy to see that $B$ is a perfect subset of $sm(S\{n\})_d$ containing $n$, and $|B|=|B_1| \leq C_{n-1}^{d-1}/2$.

Finally, by Lemma \ref{without add containing}, $A \cup B$ is a perfect subset of $sm(S)_d$, and $|A \cup B| =|A|+|B| \leq C_{n-1}^d/2 + C_{n-1}^{d-1}/2 = C_n^d/2$. This completes the proof.
\end{prof}

\vs{3mm}By Proposition \ref{homogeneous $f$-Ideal} and Theorem \ref{existence of (n, d)^{th} perfect sets}, the following corollary is clear.

\begin{cor}\label{existence of (n, d)^{th} f-ideals} For any integer $d \geq 2$ and any integer $n \geq d+2$, $V(n, d) \neq \emptyset$ if and only if $2 \mid C_n^d$.
\end{cor}

\section{ Algorithms for constructing examples of $(n, d)^{th}$ $f$-ideals }

In this section, we will show some algorithms to construct $(n, d)^{th}$ $f$-ideals. We discuss the following cases:

{\bf Case 1}: $d=2$. An $(n, 2)^{th}$ $f$-ideal is easy to construct by \cite{GWL}. For readers convenience, we repeat it as the  following: Decompose the set $[n]$ into a disjoint union of two subsets $B$ and $C$ {\it uniformly}, namely, $||B|-|C|| \leq 1$. Then set $A= \{x_ix_j \,|\, i, j \in B,$ or $i, j \in C\}$ to obtain an $(n, 2)^{th}$ perfect set. Note that $|A|=N_{(n, 2)} \leq C_n^2/2$, choose a subset $D$ of $sm(S)_2 \setminus A$ randomly, such that $|D| = C_n^2/2 - N_{(n, 2)}$ holds. It is easy to see that $A \cup D$ is still a perfect set, and $|A \cup D|=C_n^2/2$. By Proposition \ref{homogeneous $f$-Ideal}, the ideal generated by $A \cup D$ is an $(n, 2)^{th}$ $f$-ideal. Note that each $(n, 2)^{th}$ $f$-ideal can be obtained in this way except $C_5$ by \cite{GWL}.

{\bf Case 2}: $d>2$ and $n=d+2$.

\begin{alg}\label{V(d+2, d)} In order to build an $f$-ideal $I \in V(d+2, d)$, we obey the following steps:

Step 1: Calculate $C_{d+2}^d/2$. Note that $C_{d+2}^d/2=C_{d+2}^2/2$.

Step 2: As in the case 1, find a perfect subset $B$ of $sm(S)_2$ such that $|B| \leq C_{d+2}^2/2$, where $S=K[x_1, \ldots, x_{d+2}]$.

Step 3: Let $A=B^{'}$. Then $A$ is a perfect subset of $sm(S)_d$ by Lemma \ref{complement lemma}, and $|A|=|B| \leq C_{d+2}^2/2 = C_{d+2}^d/2$.

Step 4: Choose a subset $D$ of $sm(S)_d \setminus A$ randomly, such that $|D| = C_{d+2}^d/2 - |A|$ holds. It is easy to see that $M= A \cup D$ is still a perfect set, and $|A \cup D|=C_{d+2}^d/2$.

Step 5: Let $I$ be the ideal generated by $A \cup D$, By Proposition \ref{homogeneous $f$-Ideal} again, $I$ is an $(d+2, d)^{th}$ $f$-ideal.
\end{alg}

Note that in this way, we constructed almost all $(d+2, d)^{th}$ $f$-ideals.


\begin{ex}\label{V(8, 6)} {\it Show an $f$-ideal $I \in V(8, 6)$.}

Note that $8=6+2$, we obey the Algorithm \ref{V(d+2, d)}.

Note that $C_8^6/2 = 14$.
Find a perfect subset $B$ of $sm(S)_2$ such that $|B| \leq C_8^2/2 =14$, where $S=K[x_1, \ldots, x_{8}]$. It is easy to see that
$$B= \{x_1x_2, x_1x_3, x_1x_4, x_2x_3, x_2x_4, x_3x_4, x_5x_6, x_5x_7, x_5x_8, x_6x_7, x_6x_8, x_7x_8\}$$
is a perfect subset of $sm(S)_2$, with $|B|=12$. Let
$$A = B^{'}= \{x_3x_4x_5x_6x_7x_8, x_2x_4x_5x_6x_7x_8, x_2x_3x_5x_6x_7x_8, x_1x_4x_5x_6x_7x_8,$$
$$x_1x_3x_5x_6x_7x_8, x_1x_2x_5x_6x_7x_8, x_1x_2x_3x_4x_7x_8, x_1x_2x_3x_4x_6x_8,$$
$$x_1x_2x_3x_4x_6x_7, x_1x_2x_3x_4x_5x_8, x_1x_2x_3x_4x_5x_7, x_1x_2x_3x_4x_5x_6\}.$$
$A$ is a perfect subset of $sm(S)_6$. Choose $D=\{x_1x_2x_3x_5x_6x_7, x_1x_2x_4x_5x_6x_8\}$, then the ideal $I$ generated by $A \cup D$ is an $(8, 6)^{th}$ $f$-ideal.
\end{ex}

{\bf Case 3}: $d>2$ and $n > d+2$. Let $S^{[k]} = K[x_1, \ldots, x_k]$, and let $S=S^{[n]} = K[x_1, \ldots, x_n]$.

\begin{alg}\label{V(n, d)} For an integer $n > d+2$, we construct an $(n, d)^{th}$ $f$-ideal by using  the following steps:

Step 1: Let $t=n, \, l=d$ and $E = \emptyset$. Set $\mathcal{B}= \{B(t, l, E)\}$.

Step 2: Assign $\mathcal{C} = \mathcal{B}$, and denote $i = |\mathcal{C}|$.

Step 3: Choose each $B(t, l, E) \in \mathcal{C}$ one by one, deal with each one obeying the following rules:

If $l=2$ or $t = l+2$, don't change anything.

If $l \neq 2$ and $t > l+2$, then cancel $B(t, l, E)$ from $\mathcal{B}$, and add $B(t-1, l, E)$ and $B(t-1, l-1, E \cup \{t\})$ into $\mathcal{B}$.

After $i$ times, i.e., when $B(t, l, E)$ goes through all the element of $\mathcal{C}$, make a judgement:

If $l=2$ or $t=l+2$ for each $B(t, l, E) \in \mathcal{B}$, then go to step 4, else return to step 2.

Step 4: Choose $B(t, l, E) \in \mathcal{B}$ one by one, deal with each one obeying the following rules:

If $l=2$, assign $B(t, l, E)$  a perfect subset of $sm(S^{[t]})_l$ as case 1.

If $l \neq 2$ and $t = l+2$, assign $B(t, l, E)$ a perfect subset of $sm(S^{[t]})_l$ as case 2.

Step 5: For each $B(t, l, E) \in \mathcal{B}$, denote $B^{*}(t, l, E) = \{ux_{E} \,|\, u \in B(t, l, E)\}$, where $x_{E}= \prod_{j \in E} x_j$. Denote $\mathcal{B}^{*}=\cup_{B(t, l, E) \in \mathcal{B}}B^{*}(t, l, E)$.
It is direct to check that $\mathcal{B}^{*}$ is a perfect subset of $sm(S)_d$, and $|\mathcal{B}^{*}| \leq C_n^d/2$. Choose a subset $D$ of $sm(S)_d \setminus \mathcal{B}^{*}$ randomly, such that $|D| = C_n^d/2 - |\mathcal{B}^{*}|$ holds.

Step 6: Let $I$ be the ideal generated by $\mathcal{B}^{*} \cup D$.
By Proposition \ref{homogeneous $f$-Ideal} again, $I$ is an $(n, d)^{th}$ $f$-ideal.





\end{alg}

\begin{ex}\label{V(6, 3)} {\it Show a $(6, 3)^{th}$ $f$-ideal.}

Let $S=K[x_1, \ldots, x_6]$. By the above algorithm, we will choose a perfect subset $B(5, 3, \emptyset)$ of $sm(S^{5})_3$ and a perfect subset $B(5, 2, \{6\})$ of $sm(S^{5})_2$. Set $B(5, 3, \emptyset) = \{x_3x_4x_5, x_2x_4x_5, x_1x_4x_5, x_1x_2x_3\}$, and set $B(5, 2, \{6\}) = \{x_1x_2, x_1x_3, x_2x_3, x_4x_5\}$. Correspondingly, $B^{*}(5, 3, \emptyset) = B(5, 3, \emptyset)$ and $$B^{*}(5, 2, \{6\}) = \{x_1x_2x_6, x_1x_3x_6, x_2x_3x_6, x_4x_5x_6\}.$$ Hence $$\mathcal{B}^{*} = \{x_3x_4x_5, x_2x_4x_5, x_1x_4x_5, x_1x_2x_3, x_1x_2x_6, x_1x_3x_6, x_2x_3x_6, x_4x_5x_6\}$$ is a perfect subset of $sm(S)_3$. Note that $C_6^3/2 = 10$, and $|\mathcal{B}^{*}|=8$. Set $D = \{x_1x_2x_4, x_1x_2x_5\}$. The ideal $I$ generated by $\mathcal{B}^{*} \cup D$ is a $(6, 3)^{th}$ $f$-ideal.
\end{ex}

Note that the $(6, 3)^{th}$ $f$-ideal given in the above example is not unmixed. In fact, consider the simplicial complex $\sigma(sm(S)_3 \setminus G(I))$, and note that $\{1, 2\}$ is a non-face of $\sigma(sm(S)_3 \setminus G(I))$, which implies that $\sigma(sm(S)_3 \setminus G(I))$ is not a $3$-flag complex. So, $I$ is not unmixed by Proposition \ref{unmixed f-ideal}.

\section{ An upper bound of the perfect number $N_{(n, d)}$}

For a positive integer $k$ and a pair of positive integers $i \leq j$, denote by $Q_{[i, j]}^{k}$ the set of square-free monomials of degree $k$ in the polynomial ring $K[x_i, x_{i+1}, \ldots, x_j]$. Note that $Q_{[i, j]}^{k}= \emptyset$ holds for $i > j$. For a pair of monomial subsets $A$ and $B$, denote
by $A \bullet B = \{uv \,|\, u\in A, v\in B\}$. If $B = \emptyset$, then assume $A \bullet B = A$. The following theorem gives an upper bound of the $(n, d)^{th}$ perfect number for $n > d+2$.

\begin{thm}\label{N(n, d)} Given a integer $d>2$, and a integer $n\ge d+2$.  The following statements about the perfect number $N_{(n,d)}$ hold:

$(1)$ If $n=d+2$, then
\begin{equation}\label{eq: N_{(d+2,d)}}
N_{(n,d)} = N_{(n,2)} = \begin{cases}
k^2-k,  &  \text{if } n=2k; \\
k^2,  &  \text{if } n=2k+1.
\end{cases}
\end{equation}

$(2)$ If $n > d+2$, then
\begin{equation}N_{(n,d)} \leq \sum_{i=5}^{n-d+2} N_{(i, 2)}C_{n-i-1}^{d-3} + \sum_{j=3}^d N_{(j+2, 2)}C_{n-j-3}^{d-j},\end{equation}
where $C_0^0 = 1$.
\end{thm}

\begin{prof} By Lemma \ref{complement lemma} and the equation \ref{eq: perfect number new understanding} in the first section, (1) is clear.

In order to prove (2), it will suffice to show that there exists a perfect set with cardinality $t=\sum_{i=5}^{n-d+2} N_{(i, 2)}C_{n-i-1}^{d-3} + \sum_{j=3}^d N_{(j+2, 2)}C_{n-j-3}^{d-j}$.

Let $P_{(i, 2)}$ be an $(i, 2)^{th}$ perfect set with cardinality $N_{(i, 2)}$ for $5 \leq i \leq n-d+2$, and let $P_{(j+2, j)}$ be a $(j+2, j)^{th}$ perfect set with cardinality $N_{(j+2, j)}$ for $3 \leq j \leq d$. We claim that the set $$M = (\cup_{i=5}^{n-d+2}P_{(i, 2)} \bullet x_{i+1} \bullet Q_{[i+2, n]}^{d-3}) \cup (\cup_{j=3}^{d}P_{(j+2, j)} \bullet Q_{[j+4, n]}^{d-j})$$ is an $(n, d)^{th}$ perfect set, with cardinality $t$. It is easy to check that the cardinality of $M$ is $t$. It is only necessary to prove that $M$ is perfect.

For each $w \in sm(S)_{d+1}$, denote by $n_k(w)$ the cardinality of the set $\{x_i \,|\, i\leq k$ and $x_i \mid w\}$. If $n_5(w) \geq 4$, then choose the smallest $k$ such that $n_{k+3}(w)=n_{k+2}(w)=k+1$. Clearly, $3 \leq k \leq d$. It is direct to check that $w$ is divided by some monomial in $P_{(k+2, k)} \bullet Q_{[k+4, n]}^{d-k}$. If $n_5(w) \leq 3$, then choose the smallest $k$ such that $n_k(w)=3$ and $n_{k+1}(w)=4$. Clearly, $5 \leq k \leq n-d+2$. It is not hard to check that $w$ is divided by some monomial in $P_{(k, 2)} \bullet x_{k+1} \bullet Q_{[k+2, n]}^{d-3}$. Hence $M$ is upper perfect.

For each $w \in sm(S)_{d-1}$, if $n_5(w) \geq 2$, then choose the smallest $k$ such that $n_{k+3}(w)=n_{k+2}(w)=k-1$. Clearly, $3 \leq k \leq d$. It is direct to check that $w$ is covered by some monomial in $P_{(k+2, k)} \bullet Q_{[k+4, n]}^{d-k}$. If $n_5(w) \leq 1$, then choose the smallest $k$ such that $n_k(w)=1$ and $n_{k+1}(w)=2$. Clearly, $5 \leq k \leq n-d+2$ holds. It is not hard to check that $w$ is covered by some monomial in $P_{(k, 2)} \bullet x_{k+1} \bullet Q_{[k+2, n]}^{d-3}$. Hence $M$ is lower perfect.
\end{prof}

\vs{3mm}The following figure
may help to interpret the above theorem intuitively.

\begin{center}
\setlength{\unitlength}{0.11cm}
\begin{picture}(82,84)
\thinlines
\drawvector{10.0}{8.0}{64.0}{1}{0}
\drawvector{10.0}{8.0}{70.0}{0}{1}
\drawcenteredtext{6.0}{76.0}{$t$}
\drawcenteredtext{74.0}{6.0}{$l$}
\drawpath{18.0}{8.0}{18.0}{76.0}
\drawpath{10.0}{16.0}{70.0}{76.0}
\drawdotline{22.0}{76.0}{22.0}{32.0}
\drawdotline{22.0}{32.0}{64.0}{74.0}
\drawhollowdot{42.0}{76.0}
\drawdotline{42.0}{76.0}{18.0}{52.0}
\drawdotline{42.0}{76.0}{42.0}{48.0}
\drawvector{42.0}{76.0}{4.0}{0}{-1}
\drawvector{42.0}{72.0}{4.0}{-1}{-1}
\drawvector{38.0}{68.0}{4.0}{0}{-1}
\drawvector{38.0}{64.0}{4.0}{0}{-1}
\drawvector{38.0}{60.0}{4.0}{-1}{-1}
\drawvector{34.0}{56.0}{4.0}{0}{-1}
\drawvector{34.0}{52.0}{4.0}{-1}{-1}
\drawvector{30.0}{48.0}{4.0}{0}{-1}
\drawvector{30.0}{44.0}{4.0}{-1}{-1}
\drawvector{26.0}{40.0}{4.0}{0}{-1}
\drawvector{26.0}{36.0}{4.0}{0}{-1}
\drawvector{42.0}{76.0}{4.0}{-1}{-1}
\drawvector{38.0}{72.0}{4.0}{-1}{-1}
\drawvector{34.0}{68.0}{4.0}{0}{-1}
\drawvector{34.0}{64.0}{4.0}{-1}{-1}
\drawvector{30.0}{60.0}{4.0}{-1}{-1}
\drawvector{26.0}{56.0}{4.0}{0}{-1}
\drawvector{26.0}{52.0}{4.0}{-1}{-1}
\drawvector{22.0}{48.0}{4.0}{0}{-1}
\drawvector{22.0}{44.0}{4.0}{-1}{-1}
\drawcenteredtext{46.0}{78.0}{$(d, n)$}
\drawcenteredtext{66.0}{66.0}{$t=l+2$}
\drawcenteredtext{62.0}{74.0}{$t=l+3$}
\drawcenteredtext{16.0}{78.0}{$l=2$}
\drawcenteredtext{26.0}{76.0}{$l=3$}
\drawcenteredtext{8.0}{6.0}{$O$}
\drawcenteredtext{22.0}{22.0}{$(2, 4)$}
\drawcenteredtext{22.0}{32.0}{$(3, 6)$}
\end{picture}
\centerline{\fs{Figure 1. Upper Bound}}
\end{center}

In this figure, there is a boundary consisting of the line $l=2$ and the line $t=l+2$. From the point $(d, n)$ to a point of the boundary, every directed chain $\mathcal{C}$ denotes a set of monomials $M(\mathcal{C})$ by the following rules:

$(1)$ Every arrow of $\mathcal{C}$ is from $(l, t)$ to either $(l, t-1)$ or $(l-1, t-1)$.

$(2)$ If the arrow is from $(l, t)$ to $(l, t-1)$, then each monomial in $M(\mathcal{C})$ is not divided by $x_t$. Correspondingly, if it is from $(l, t)$ to $(l-1, t-1)$, then each monomial in $M(\mathcal{C})$ is divided by $x_t$.

$(3)$ Each point $(l, t)$ of the boundary is a $(t, l)^{th}$ perfect set.

Actually, the figure shows us a class of $(n, d)^{th}$ perfect sets. For each point $(l, t)$ of the boundary, if we choose the corresponding perfect set to be a $(t, l)^{th}$ perfect set with cardinality $N_{(t, l)}$, then the cardinality of the $(n, d)^{th}$ perfect set is exactly $\sum_{i=5}^{n-d+2} N_{(i, 2)}C_{n-i-1}^{d-3} + \sum_{j=3}^d N_{(j+2, 2)}C_{n-j-3}^{d-j}$.

\begin{ex}\label{N_{(6, 3)}=7} {\it Calculation of the $(6, 3)^{th}$ perfect number.}

 Let $A$ be a $(6, 3)^{th}$ perfect set.  By Proposition \ref{upper without lower containing}(2),  $A\{\check{6}\}$ is an upper perfect set without $6$. Hence $|A\{\check{6}\}| \geq N_{(5, 3)}=4$. By Proposition \ref{upper without lower containing}(3), $A\{6\}$ is a lower perfect set containing $6$. Note that for the monomials of $\{x_1, x_2, x_3, x_4, x_5\}$, each monomial in $A\{6\}$ covers at most two of them. So, $|A\{6\}| \geq 3$. Hence $|A| \geq |A\{\check{6}\}| + |A\{6\}| \geq 7$.  Actually, as showed in Example \ref{not upper containing example}, there exists a $(6, 3)^{th}$ perfect set
 $$B = \{x_1x_2x_3, x_1x_2x_4, x_1x_2x_5, x_3x_4x_5, x_1x_3x_6, x_2x_3x_6, x_4x_5x_6\}$$ with cardinality $7$. Thus $N_{(6, 3)} = 7$. Note that  the upper bound given by Proposition \ref{N(n, d)}(2) is $8$, and is not bad for  the perfect number in the case.
\end{ex}


\section{Nonhomogeneous $f$-ideal}

In \cite{GWL}, a characterization of $f$-ideals in general case is shown, but it is still not easy to show an example of nonhomogeneous $f$-ideal. In fact, the interference from monomials of different degree makes the computation complicated. Anyway, we worked out the following examples:

\begin{ex}\label{nonhomogeneous}Let $S=K[x_1, x_2, x_3, x_4, x_5]$, and let $$I = \langle x_1x_2, x_3x_4, x_1x_3x_5, x_2x_4x_5 \rangle.$$
It is direct to check that $$\delta_{\mathcal{F}}(I) = \langle \{1, 2\}, \{3, 4\}, \{1, 3, 5\}, \{2, 4, 5\} \rangle$$
and $$\delta_{\mathcal{N}}(I) = \langle \{1, 3\}, \{2, 4\}, \{1, 4, 5\}, \{2, 3, 5\} \rangle. $$
It is easy to see they have the same $f$-vector, and hence $I$ is an  $f$-ideal, which is clearly nonhomogeneous.
\end{ex}

\vs{3mm}In fact, there are a lot of nonhomogeneous $f$-ideals. We will show another example to end this section.

\begin{ex}\label{nonhomogeneous2}Let $S=K[x_1, x_2, x_3, x_4, x_5, x_6]$, and let $$I = \langle x_1x_2, x_2x_3, x_1x_3, x_4x_5, x_1x_4x_6, x_1x_5x_6, x_2x_4x_6 \rangle.$$
Note that $$\delta_{\mathcal{N}}(I) = \langle \{1, 4\}, \{1, 5\}, \{1, 6\}, \{2, 4\}, \{2, 5, 6\}, \{3, 4, 6\}, \{3, 5, 6\} \rangle. $$
It is direct to check that $I$ is also a nonhomogeneous $f$-ideal.
\end{ex}

\end{document}